\documentclass[12pt]{article}

\usepackage{amsmath,amssymb,amsthm}
\usepackage{enumitem}
\usepackage[a4paper,margin=2.5cm]{geometry}
\newtheorem{The}{Theorem}[section]
\newtheorem{Lem}[The]{Lemma}
\newtheorem{Prop}[The]{Proposition}
\newtheorem{Cor}[The]{Corollary}
\newtheorem{Question}[The]{Question}
\theoremstyle{definition}
\newtheorem{Def}[The]{Definition}
\newtheorem{Ex}[The]{Example}

\newcommand{\C}{\mathbb{C}}

\newcommand{\Z}{\mathbb{Z}}

\begin{document}

\begin{center}
{\Large\bf Hyperbolicity Notions on Complete Non-Compact Complex Manifolds}
\end{center}

\vspace{2ex}

\begin{center}
{\large Samir Marouani}
\end{center}

\vspace{3ex}

\noindent{\small{\bf Abstract.} Gromov's K\"ahler hyperbolicity, extended by McMullen to complete manifolds of finite volume and bounded geometry, has non-K\"ahler analogues (balanced, SKT, Gauduchon and divisorial hyperbolicity) that have so far been studied only on compact manifolds. We define these notions for complete Hermitian manifolds and prove the following. A manifold carrying a balanced metric and a Hermitian metric with bounded $(\partial+\bar\partial)$-potentials satisfies a Cheeger inequality, so its fundamental group is non-amenable in the compact case. K\"ahler and SKT hyperbolic manifolds are Brody hyperbolic, and balanced and Gauduchon hyperbolic manifolds are divisorially hyperbolic, without any compactness assumption. Products give complete non-compact non-K\"ahler examples. A non-compact surface may be divisorially hyperbolic without being Brody hyperbolic. Finally, no functorial pseudo-distance characterises divisorial hyperbolicity.}

\vspace{4ex}

\section{Introduction}\label{section:Introduction}

Gromov introduced in [Gro91] the notion of {\bf K\"ahler hyperbolicity}: a compact complex manifold $X$ is K\"ahler hyperbolic if it carries a K\"ahler metric $\omega$ such that the lifted form $\widetilde\omega=\pi_X^\star\omega$ to the universal cover $\pi_X:\widetilde{X}\to X$ is $d(\mbox{bounded})$. The condition is strong. Such a manifold is Kobayashi hyperbolic, its $L^2$ cohomology on $\widetilde X$ is concentrated in the middle degree, and the sign of its Euler characteristic is $(-1)^{\mbox{dim}_\C X}$. Gromov also showed [Gro91, Corollary 0.4.C] that $K_X$ is big, and Chen and Yang proved [CY17, Theorem 2.11] that $K_X$ is ample.

Gromov's definition makes sense for complete K\"ahler manifolds, provided one asks for finite volume and bounded geometry. McMullen used this in [McM00] to show that the moduli space $\mathcal{M}_{g,n}$ of Riemann surfaces is K\"ahler hyperbolic for a suitable complete K\"ahler metric. This is remarkable because, in dimension greater than one, $\mathcal{M}_{g,n}$ carries no complete metric of bounded negative curvature.

Several non-K\"ahler analogues of K\"ahler hyperbolicity have been introduced for {\it compact} manifolds: {\bf balanced hyperbolicity} and {\bf divisorial hyperbolicity} in [MP22], {\bf SKT hyperbolicity} and {\bf Gauduchon hyperbolicity} in [Mar24], and Hermitian-symplectic hyperbolicity in [Khe26]. In [Mar26] we showed in particular that divisorial and Brody hyperbolicity coincide for compact surfaces, that a compact K\"ahler manifold which is Gauduchon hyperbolic has a non-amenable fundamental group, and that a Gauduchon hyperbolic manifold contains no compact hypersurface lifting to its universal cover.

The aim of this paper is to ask what survives of these notions, and of the arguments behind their consequences, once compactness is dropped. This is precisely the setting of McMullen's theorem. We say that a Hermitian manifold $(X,\omega)$ has {\bf full bounded geometry} if it is complete, of finite volume, has bounded sectional curvature, and the injectivity radius of $(\widetilde{X},\widetilde\omega)$ is bounded below (Definition \ref{Def:bounded-geometry}). These conditions are automatic when $X$ is compact. In \S\ref{section:definitions} we define the K\"ahler, SKT, balanced and Gauduchon hyperbolic metrics on manifolds with full bounded geometry, by adding McMullen's conditions to the compact definitions.

\vspace{2ex}

\noindent {\bf Main results.} The first result is a Cheeger-type inequality. It pairs a balanced metric with a Hermitian metric whose K\"ahler form has bounded $(\partial+\bar\partial)$-potentials, as is the case for a manifold that is both K\"ahler and SKT hyperbolic.

\begin{The}\label{The:Cheeger} Let $X$ be a complex manifold of dimension $n\geq 2$ carrying a balanced metric $\omega_B$ and a Hermitian metric $\omega_S$, uniformly equivalent to $\omega_B$, such that $\widetilde\omega_S=\partial\alpha+\bar\partial\beta$ on $\widetilde{X}$ for some bounded forms $\alpha,\beta$. Then there exists a constant $C>0$ such that $$\mbox{vol}(F)\leq C\,\mbox{area}(\partial F)$$ for every relatively compact open set $F\subset\widetilde{X}$ with smooth boundary, the volume and the area being taken for $\widetilde g_B$ (and, with another constant, for $\widetilde g_S$). Consequently, the Cheeger constant of $\widetilde{X}$ is positive and $\lambda_0(\widetilde{X})>0$. If $X$ is compact, $\pi_1(X)$ is non-amenable.\end{The}

In particular, a compact manifold carrying a balanced metric and having an amenable fundamental group is never SKT hyperbolic. This applies, for example, to $\#_k(S^3\times S^3)$ and to the smooth ample hypersurfaces of simple abelian varieties (\S\ref{section:Cheeger}).

The second result says that the known implications to Brody and divisorial hyperbolicity do not require compactness.

\begin{The}\label{The:Brody-div} Let $(X,\omega)$ be a Hermitian manifold of dimension $n\geq 2$, with no further assumption.

\vspace{1ex}

(i)\, If $\omega$ is $\widetilde{d}(\mbox{bounded})$ or $\widetilde{\partial+\bar\partial}(\mbox{bounded})$, every holomorphic map $\C\to X$ is constant.

\vspace{1ex}

(ii)\, If $\omega^{n-1}$ is $\widetilde{d}(\mbox{bounded})$ or $\widetilde{\partial+\bar\partial}(\mbox{bounded})$, there is no holomorphic map $f:\C^{n-1}\to X$ that is non-degenerate at some point and of subexponential growth with respect to $\omega$ in the sense of [MP22, Definition 2.3], and no holomorphic map $\mathbb{P}^{n-1}\to X$ that is non-degenerate at some point.\end{The}

Thus K\"ahler hyperbolic and SKT hyperbolic manifolds are Brody hyperbolic, and balanced hyperbolic and Gauduchon hyperbolic manifolds are divisorially hyperbolic, whether compact or not. In the compact case these are [Gro91], [Mar24, Theorem 3.5], [MP22, Theorem 2.8] and [Mar24, Theorem 3.11]. Next, products give examples that are complete, of finite volume and non-K\"ahler.

\begin{The}\label{The:products} If $(X_1,\omega_1)$ and $(X_2,\omega_2)$ are balanced hyperbolic of dimensions $n_1,n_2\geq 2$, then $(X_1\times X_2,\,\omega_1+\omega_2)$ is balanced hyperbolic. The same holds for Gauduchon hyperbolicity. In particular, let $X_1=\mathcal{M}_{0,5}$ be the moduli space of five labelled points of $\mathbb{P}^1$ and $X_2=\#_2(S^3\times S^3)$ with a Friedman--Lu--Tian complex structure. Then $X_1\times X_2$ has full bounded geometry, is balanced hyperbolic, carries no K\"ahler metric, and has $\lambda_0>0$ on its universal cover.\end{The}

For compact surfaces, divisorial and Brody hyperbolicity coincide [Mar26, Proposition 3.1]. Compactness is really needed for this.

\begin{The}\label{The:surface} The complete K\"ahler surface $X=\C\times\mathbb{D}$, with the metric $\omega=e^{2|z|^2}|dz|^2+\omega_{\mathbb{D}}$ ($\omega_{\mathbb{D}}$ the Poincar\'e metric), is divisorially hyperbolic with respect to $\omega$ but is not Brody hyperbolic. It is not divisorially hyperbolic with respect to $|dz|^2+\omega_{\mathbb{D}}$.\end{The}

Finally, one may ask for an analogue of the Kobayashi pseudo-distance adapted to divisorial hyperbolicity.

\begin{The}\label{The:pseudo-distance} Let $X\mapsto d_X$ be a rule assigning to every complex manifold a continuous pseudo-distance such that holomorphic maps do not increase distances. Then $d_\C\equiv 0$, and if $d_X$ is a distance, $X$ is Brody hyperbolic. In particular, there is no such rule for which $d_X$ is a distance if and only if $X$ is divisorially hyperbolic, even for compact $X$.\end{The}

\vspace{2ex}

\noindent {\bf Complements.} A few further statements are proved in \S\ref{section:complements}. For $1\leq p\leq n-1$, a $(p,\,p)$-form $\Omega$ on $(X,\omega)$ is said to be {\it uniformly positive} if it is bounded and $\Omega|_P\geq c\,\omega^p|_P$ on every complex $p$-dimensional subspace $P$ of a tangent space, for some $c>0$. We say that $(X,\omega,\Omega)$ is {\bf $p$-SKT hyperbolic} if $\widetilde\Omega=\partial\alpha+\bar\partial\beta$ with $\alpha,\beta$ bounded; $p=1$, $\Omega=\omega$ is SKT hyperbolicity and $p=n-1$, $\Omega=\omega^{n-1}$ is Gauduchon hyperbolicity.

\vspace{1ex}

(a)\, (Proposition \ref{Prop:isoperimetric}) A $p$-SKT hyperbolic manifold satisfies $\mbox{vol}_{2p}(Z)\leq C\,\mbox{vol}_{2p-1}(\partial Z)$ for the complex $p$-dimensional submanifolds $Z$ of $\widetilde{X}$.

\vspace{1ex}

(b)\, (Corollary \ref{Cor:obstruction}) If a compact manifold $X$ contains a compact complex $p$-dimensional submanifold $Z$ such that the image of $\pi_1(Z)\to\pi_1(X)$ is amenable, then $X$ is not $p$-SKT hyperbolic. For $p=n-1$ this extends [Mar26, Proposition 5.1], which treats the hypersurfaces with trivial image.

\vspace{1ex}

(c)\, (Proposition \ref{Prop:Poincare}) If $\widetilde\omega^{n-1}=d\beta$ with $\beta$ bounded, then $\|df\|_{L^2}\leq C\,\|dd^cf\|_{L^2}$ for all $f\in C^\infty_c(\widetilde{X})$. Moreover, $\lambda_0(\widetilde{X})>0$ under a smallness condition on the torsion term $d\widetilde\omega\wedge\beta$.

\vspace{2ex}

We also say clearly what is not done. We do not construct complete finite-volume balanced hyperbolic manifolds that are neither K\"ahler hyperbolic nor products; McMullen's construction, which corrects an incomplete form by $\partial\bar\partial\,\mbox{Log}(\varepsilon/d)$, would require the control of $d(\partial\bar\partial h\wedge\omega^{n-2})$, which vanishes when $\omega$ is K\"ahler but not in general. We do not know whether K\"ahler hyperbolicity in the complete sense implies Kobayashi hyperbolicity, since Theorem \ref{The:Brody-div} only gives Brody hyperbolicity (see Question \ref{Question:Kobayashi}), nor whether Brody and divisorial hyperbolicity coincide for non-compact surfaces of finite volume and bounded geometry (Question \ref{Question:surfaces}).

\section{Hyperbolicity notions on complete manifolds}\label{section:definitions}

Throughout the paper, $X$ is a connected complex manifold of dimension $n\geq 2$, {\it not necessarily compact}, $\omega$ is a Hermitian metric on $X$ (a $C^\infty$ positive definite $(1,\,1)$-form) and $g$ is the associated Riemannian metric. We denote by $\pi_X:\widetilde{X}\to X$ the universal cover and put $\widetilde\omega=\pi_X^\star\omega$ and $\widetilde g=\pi_X^\star g$. On $k$-forms we use the norm induced by the metric (the Hilbert--Schmidt norm), so that $|\varphi(w_1,\dots,w_k)|\leq|\varphi|\,|w_1|\cdots|w_k|$. A form on $\widetilde{X}$ is {\it bounded} if its sup-norm for $\widetilde\omega$ is finite.

\begin{Def}\label{Def:bounded-geometry} The Hermitian manifold $(X,\omega)$ has {\bf full bounded geometry} if:

\vspace{1ex}

(G1)\, $g$ is complete;

(G2)\, $\mbox{vol}(X,g)<\infty$;

(G3)\, the sectional curvature of $g$ is bounded above and below;

(G4)\, the injectivity radius of $(\widetilde{X},\widetilde g)$ is bounded below by a positive constant.

\vspace{1ex}

\noindent These conditions hold automatically if $X$ is compact.\end{Def}

Recall that a $C^\infty$ $k$-form $\varphi$ on $X$ is said to be $\widetilde{d}(\mbox{bounded})$ (w.r.t. $\omega$) if $\pi_X^\star\varphi=d\beta$ for some bounded $(k-1)$-form $\beta$ on $\widetilde{X}$ [Gro91]. We shall also use the following notion.

\begin{Def}\label{Def:bounded-potentials} A $C^\infty$ $k$-form $\varphi$ on $X$ is said to be $\widetilde{\partial+\bar\partial}(\mbox{bounded})$ (w.r.t. $\omega$) if $\pi_X^\star\varphi=\partial\alpha+\bar\partial\beta$ for some bounded $(k-1)$-forms $\alpha,\beta$ on $\widetilde{X}$ [Mar24].\end{Def}

\begin{Def}\label{Def:hyperbolic-notions} Let $(X,\omega)$ have full bounded geometry. The metric $\omega$ is said to be:

\vspace{1ex}

$\bullet$\, {\bf K\"ahler hyperbolic} if $d\omega=0$ and $\omega$ is $\widetilde{d}(\mbox{bounded})$;

$\bullet$\, {\bf SKT hyperbolic} if $\partial\bar\partial\omega=0$ and $\omega$ is $\widetilde{\partial+\bar\partial}(\mbox{bounded})$;

$\bullet$\, {\bf balanced hyperbolic} if $d\omega^{n-1}=0$ and $\omega^{n-1}$ is $\widetilde{d}(\mbox{bounded})$;

$\bullet$\, {\bf Gauduchon hyperbolic} if $\partial\bar\partial\omega^{n-1}=0$ and $\omega^{n-1}$ is $\widetilde{\partial+\bar\partial}(\mbox{bounded})$;

$\bullet$\, {\bf degenerate balanced (bounded)} if $d\omega^{n-1}=0$ and $\omega^{n-1}=d\beta$ on $X$ itself for some bounded $(2n-3)$-form $\beta$.

\vspace{1ex}

\noindent The manifold $X$ is said to be of the corresponding type if it carries such a metric.\end{Def}

For compact $X$ these are the definitions of [Gro91], [Mar24] and [MP22]; the last one is the class of degenerate balanced metrics of [MP22, Definition 2.11]. Since boundedness is measured with respect to $\omega$, these notions depend on the metric and not only on the complex structure, in contrast with the compact case.

\begin{Lem}\label{Lem:elementary} (a)\, If $\varphi=\sum\varphi^{p,q}$ is the decomposition of a form into bidegrees, then $|\varphi^{p,q}|\leq|\varphi|$. In particular, the components of a bounded form are bounded.

\vspace{1ex}

(b)\, If $c^{-1}g_1\leq g_2\leq c\,g_1$ for two Riemannian metrics, the pointwise norms of $k$-forms for $g_1$ and $g_2$ differ by a factor at most $c^{k/2}$, and the $m$-dimensional volumes by a factor at most $c^{m/2}$.

\vspace{1ex}

(c)\, Let $\varphi_j$ be a form on $N_j$ ($j=1,2$) and let $N_1\times N_2$ carry the product metric. Then $|\varphi_1\wedge\varphi_2|=|\varphi_1|\,|\varphi_2|$ for the pull-backs to the product.

\vspace{1ex}

(d)\, For forms $\varphi,\psi$ of degrees $k,l$ on a Riemannian manifold of dimension $m$, $|\varphi\wedge\psi|\leq C(m,k,l)\,|\varphi|\,|\psi|$.\end{Lem}

\noindent {\it Proof.} (a) The splitting by bidegree is orthogonal for a Hermitian metric. (b) Diagonalise $g_2$ with respect to $g_1$ at a point. The eigenvalues lie in $[c^{-1},c]$, the squared norm of a basis $k$-form $dx^I$ is multiplied by a product of $k$ eigenvalues (or of their inverses), and the volume form by the square root of the product of all eigenvalues. (c) $\Lambda^k(V_1\oplus V_2)^\star$ is the orthogonal direct sum of the spaces $\Lambda^aV_1^\star\otimes\Lambda^bV_2^\star$, and the norm is multiplicative on pure tensors. (d) In an orthonormal coframe, the coefficients of $\varphi\wedge\psi$ are signed sums of products $\varphi_I\psi_J$ over the splittings of a multi-index; by the Cauchy--Schwarz inequality, each coefficient is bounded by a combinatorial constant times $|\varphi|\,|\psi|$. \hfill $\qed$

\begin{Prop}\label{Prop:implications} Let $n\geq 2$ and keep the same metric $\omega$. If $\omega$ is K\"ahler hyperbolic, it is SKT hyperbolic and balanced hyperbolic. If $\omega$ is balanced hyperbolic, it is Gauduchon hyperbolic. If $\omega$ is degenerate balanced (bounded), it is balanced hyperbolic.\end{Prop}

\noindent {\it Proof.} Let $\omega$ be K\"ahler hyperbolic, with $\widetilde\omega=d\theta$ and $\theta=\theta^{1,0}+\theta^{0,1}$ bounded (Lemma \ref{Lem:elementary}(a)). Expanding by bidegree, $d\theta=\partial\theta^{1,0}+(\bar\partial\theta^{1,0}+\partial\theta^{0,1})+\bar\partial\theta^{0,1}$. Since $\widetilde\omega$ has pure type $(1,\,1)$, we get $\widetilde\omega=\partial\theta^{0,1}+\bar\partial\theta^{1,0}$ and $\partial\bar\partial\omega=0$: $\omega$ is SKT hyperbolic. Moreover, $d\omega^{n-1}=(n-1)\omega^{n-2}\wedge d\omega=0$ and $\widetilde\omega^{n-1}=d(\theta\wedge\widetilde\omega^{n-2})$. The primitive is bounded because $|\widetilde\omega|$ is a constant depending only on $n$ and $|\theta\wedge\widetilde\omega^{n-2}|\leq C|\theta|\,|\widetilde\omega|^{n-2}$ (Lemma \ref{Lem:elementary}(d)). So $\omega$ is balanced hyperbolic (this is [MP22, Lemma 2.2]).

Let $\omega$ be balanced hyperbolic, with $\widetilde\omega^{n-1}=d\beta=\partial\beta+\bar\partial\beta$. The potentials are bounded, and $d\omega^{n-1}=0$ forces $\partial\omega^{n-1}=\bar\partial\omega^{n-1}=0$ (the two parts have different types), hence $\partial\bar\partial\omega^{n-1}=0$.

If $\omega^{n-1}=d\beta$ on $X$ with $\beta$ bounded, then $\widetilde\omega^{n-1}=d(\pi_X^\star\beta)$ and $\pi_X^\star\beta$ is bounded because $\pi_X$ is a local isometry. \hfill $\qed$

\section{A Cheeger inequality and non-amenability}\label{section:Cheeger}

For a Riemannian manifold $N$, put $$h(N)=\inf_F\frac{\mbox{area}(\partial F)}{\mbox{vol}(F)},\hspace{4ex}\lambda_0(N)=\inf_{f\in C^\infty_c(N)\setminus\{0\}}\frac{\int_N|df|^2}{\int_Nf^2},$$ where $F$ runs over the relatively compact open subsets of $N$ with smooth boundary. Cheeger's inequality [Che70] states that $\lambda_0(N)\geq h(N)^2/4$, with no hypothesis on $N$. We shall also use Brooks' theorem [Bro81]: if $X$ is compact, then $\lambda_0(\widetilde{X})>0$ if and only if $\pi_1(X)$ is non-amenable.

The proof of Theorem \ref{The:Cheeger} rests on the following pointwise lemma.

\begin{Lem}\label{Lem:bidegree} Let $\Phi$ be a form of bidegree $(n-1,\,n-1)$ on a complex manifold of dimension $n$ and let $\xi$ be a $1$-form. Then $$\partial\xi\wedge\Phi=\partial\xi^{0,1}\wedge\Phi\hspace{3ex}\mbox{and}\hspace{3ex}\bar\partial\xi\wedge\Phi=\bar\partial\xi^{1,0}\wedge\Phi.$$ If, moreover, $d\Phi=0$, then $\partial\xi^{0,1}\wedge\Phi=d(\xi^{0,1}\wedge\Phi)$ and $\bar\partial\xi^{1,0}\wedge\Phi=d(\xi^{1,0}\wedge\Phi)$.\end{Lem}

\noindent {\it Proof.} The form $\partial\xi^{1,0}$ has type $(2,0)$, and $(2,0)+(n-1,n-1)=(n+1,n-1)$ is not a possible bidegree, so $\partial\xi^{1,0}\wedge\Phi=0$; this gives the first identity, and the second is symmetric. Next, $\bar\partial\xi^{0,1}\wedge\Phi$ has bidegree $(n-1,n+1)$, hence vanishes, so $\partial\xi^{0,1}\wedge\Phi=d\xi^{0,1}\wedge\Phi$. Finally, $d(\xi^{0,1}\wedge\Phi)=d\xi^{0,1}\wedge\Phi-\xi^{0,1}\wedge d\Phi=d\xi^{0,1}\wedge\Phi$ if $d\Phi=0$. \hfill $\qed$

\vspace{2ex}

\noindent {\it Proof of Theorem \ref{The:Cheeger}.} Write $\omega_B$ for $\widetilde\omega_B$, put $\Theta=\widetilde\omega_S=\partial\alpha+\bar\partial\beta$, and fix $c\geq 1$ with $c^{-1}\omega_B\leq\omega_S\leq c\,\omega_B$. Let $F\subset\widetilde{X}$ be a relatively compact open set with smooth boundary. We estimate $\int_F\Theta\wedge\omega_B^{n-1}$ from below and from above.

\vspace{1ex}

$\bullet$\, {\it Lower bound.} For a real $(1,\,1)$-form $\lambda$ one has $\lambda\wedge\omega_B^{n-1}=\frac{1}{n}(\mbox{tr}_{\omega_B}\lambda)\,\omega_B^n$. Since $\Theta\geq c^{-1}\omega_B$, its trace is at least $n\,c^{-1}$, hence $\Theta\wedge\omega_B^{n-1}\geq c^{-1}\omega_B^n=c^{-1}n!\,dV_B$ and $$\int_F\Theta\wedge\omega_B^{n-1}\geq c^{-1}\,n!\,\mbox{vol}_B(F).$$

\vspace{1ex}

$\bullet$\, {\it Upper bound.} Lemma \ref{Lem:bidegree} applied with $\Phi=\omega_B^{n-1}$, which is closed because $\omega_B$ is balanced, gives $$\Theta\wedge\omega_B^{n-1}=\partial\alpha\wedge\omega_B^{n-1}+\bar\partial\beta\wedge\omega_B^{n-1}=d\big((\alpha^{0,1}+\beta^{1,0})\wedge\omega_B^{n-1}\big).$$ By Stokes' theorem, $$\int_F\Theta\wedge\omega_B^{n-1}=\int_{\partial F}(\alpha^{0,1}+\beta^{1,0})\wedge\omega_B^{n-1}\leq C'\,\mbox{area}_B(\partial F).$$ Here $C'$ depends on the sup-norms of $\alpha,\beta$ (which are bounded for $\widetilde g_S$, hence for $\widetilde g_B$ by Lemma \ref{Lem:elementary}(b), and so are their components by Lemma \ref{Lem:elementary}(a)) and on the constant norm of $\omega_B^{n-1}$, together with Lemma \ref{Lem:elementary}(d).

\vspace{1ex}

Comparing the two estimates yields $\mbox{vol}_B(F)\leq C\,\mbox{area}_B(\partial F)$. Hence $h(\widetilde{X})\geq 1/C$ and $\lambda_0(\widetilde{X})\geq 1/(4C^2)$ by Cheeger's inequality. If $X$ is compact, any two metrics on $X$ are uniformly equivalent, and Brooks' theorem shows that $\pi_1(X)$ is non-amenable. \hfill $\qed$

\vspace{2ex}

\begin{Cor}\label{Cor:McMullen} Let $(X,\omega)$ be K\"ahler hyperbolic. Then $h(\widetilde{X})>0$ and $\lambda_0(\widetilde{X})>0$ for the metric $\widetilde g$.\end{Cor}

\noindent {\it Proof.} Apply Theorem \ref{The:Cheeger} with $\omega_S=\omega_B=\omega$. A K\"ahler metric is balanced, and by the proof of Proposition \ref{Prop:implications}, $\widetilde\omega=\partial\theta^{0,1}+\bar\partial\theta^{1,0}$ with $\theta$ bounded, hence with bounded components. \hfill $\qed$

\vspace{2ex}

The condition $\partial\bar\partial\omega_S=0$ was never used in the proof of Theorem \ref{The:Cheeger}, and neither was completeness. Lemma 5.2 of [Mar26] uses a different pairing: a K\"ahler metric against the potentials of $\omega^{n-1}$ for a compact Gauduchon hyperbolic $\omega$. Theorem \ref{The:Cheeger} pairs a balanced metric with the potentials of $\omega_S$ itself, and neither statement contains the other.

\begin{Ex}\label{Ex:not-SKT-hyperbolic} (a)\, {\it Ample hypersurfaces in simple abelian varieties.} Let $A$ be a simple abelian variety of dimension $n\geq 3$ and let $X\subset A$ be a smooth ample hypersurface (it exists by Bertini's theorem). We check four facts.

\vspace{1ex}

(1)\, $X$ is Brody hyperbolic. By the theorem of Bloch--Ochiai--Green [Och77], [Gre75], the Zariski closure of a non-constant entire curve in $X$ is a translate of a positive-dimensional abelian subvariety of $A$ contained in $X$. As $A$ is simple, such a subvariety is $A$ itself, which is not contained in $X$.

\vspace{1ex}

(2)\, $\pi_1(X)\cong\Z^{2n}$. Indeed, $X$ is ample with $\mbox{dim}_\C X=n-1\geq 2$, so the Lefschetz hyperplane theorem gives $\pi_1(X)\cong\pi_1(A)=\Z^{2n}$, which is amenable.

\vspace{1ex}

(3)\, $X$ carries a balanced metric: it is projective, so the restriction of a Fubini--Study form is a K\"ahler form, and a K\"ahler form is balanced.

\vspace{1ex}

(4)\, If $X$ were SKT hyperbolic, with an SKT hyperbolic metric $\omega_S$, then Theorem \ref{The:Cheeger}, applied to a K\"ahler metric $\omega_B$ and to $\omega_S$ (which are uniformly equivalent because $X$ is compact), would make $\pi_1(X)$ non-amenable, contradicting (2).

\vspace{1ex}

\noindent A concrete instance is $A=J(C)$ for a very general curve $C$ of genus $3$: $J(C)$ is simple, $C$ is not hyperelliptic, and the theta divisor $\Theta=W_2\cong C^{(2)}$ is a smooth ample surface with $\pi_1=\Z^6$.

\vspace{2ex}

(b)\, {\it The manifolds $X_k$.} Let $k\geq 2$ and let $X_k=\#_k(S^3\times S^3)$ with a Friedman--Lu--Tian complex structure [Fri89], [LT93].

\vspace{1ex}

(1)\, $X_k$ is simply connected (Lemma \ref{Lem:X_k}) and carries a balanced metric $\omega_B$ [FLY12], see also [MP22, Theorem 1.2].

\vspace{1ex}

(2)\, Suppose $\omega_S$ is a Hermitian metric on $X_k$ with $\omega_S=\partial\alpha+\bar\partial\beta$. Since $\widetilde{X}_k=X_k$ is compact, $\alpha$ and $\beta$ are automatically bounded.

\vspace{1ex}

(3)\, Lemma \ref{Lem:bidegree}, applied with $n=3$ and $\Phi=\omega_B^2$, gives $\omega_S\wedge\omega_B^2=d\big((\alpha^{0,1}+\beta^{1,0})\wedge\omega_B^2\big)$. Hence $\int_{X_k}\omega_S\wedge\omega_B^2=0$. But $\omega_S\wedge\omega_B^2$ is a positive volume form, a contradiction.

\vspace{1ex}

\noindent So $X_k$ is not SKT hyperbolic, although it is balanced hyperbolic (Example \ref{Ex:product}).\end{Ex}

The second example shows that the amenability obstruction, which comes from the closedness of the K\"ahler form, has no counterpart for balanced hyperbolicity.

\begin{Question}\label{Question:Kaehler-SKT} Does a compact manifold carrying both a K\"ahler metric and an SKT hyperbolic metric admit a K\"ahler hyperbolic metric? By Theorem \ref{The:Cheeger}, its fundamental group is non-amenable.\end{Question}

\section{Products and non-K\"ahler examples}\label{section:products}

\begin{Prop}\label{Prop:volume-form} Let $(N,g)$ be an oriented Riemannian manifold of dimension $m$ whose volume form satisfies $dV=d\eta$ with $\eta$ bounded. Then $h(N)\geq 1/\|\eta\|_\infty$ and $\lambda_0(N)\geq 1/(4\|\eta\|_\infty^2)$.\end{Prop}

\noindent {\it Proof.} For a relatively compact open set $F$ with smooth boundary, $\mbox{vol}(F)=\int_{\partial F}\eta\leq\|\eta\|_\infty\,\mbox{area}(\partial F)$. Apply Cheeger's inequality. \hfill $\qed$

\begin{Lem}\label{Lem:product-geometry} Let $X_1,X_2$ be complete Riemannian manifolds with bounded sectional curvature, finite volume, and injectivity radius of the universal cover bounded below. Then the same holds for $X_1\times X_2$ with the product metric.\end{Lem}

\noindent {\it Proof.} $\bullet$\, {\it Completeness and volume.} The product distance is $\sqrt{d_1^2+d_2^2}$, so Cauchy sequences project to Cauchy sequences in each factor and $X_1\times X_2$ is complete. Moreover, $\mbox{vol}(X_1\times X_2)=\mbox{vol}(X_1)\,\mbox{vol}(X_2)<\infty$.

\vspace{1ex}

$\bullet$\, {\it Curvature.} For tangent vectors $u=u_1+u_2$ and $v=v_1+v_2$, the curvature tensor of the product satisfies $R(u,v,v,u)=R_1(u_1,v_1,v_1,u_1)+R_2(u_2,v_2,v_2,u_2)$. If $K_j$ denotes the sectional curvature of $X_j$ on the plane spanned by $u_j,v_j$ (the term is zero when $u_j,v_j$ are dependent), then $$\sec(u,v)\,|u\wedge v|^2=K_1|u_1\wedge v_1|^2+K_2|u_2\wedge v_2|^2.$$ Since $u\wedge v$ decomposes orthogonally along $\Lambda^2V_1$, $\Lambda^2V_2$ and $V_1\otimes V_2$, we have $|u_1\wedge v_1|^2+|u_2\wedge v_2|^2\leq|u\wedge v|^2$. Hence $\sec(u,v)$ lies between $\min(0,\inf K_j)$ and $\max(0,\sup K_j)$, which proves (G3).

\vspace{1ex}

$\bullet$\, {\it Injectivity radius.} The universal cover of the product is the product of the universal covers. A geodesic of a product is a pair of geodesics parametrised proportionally to arc length, and it stops minimising exactly when one of its components does; so the injectivity radius of the product is the minimum of those of the factors, which gives (G4). \hfill $\qed$

\begin{The}\label{The:product-theorem} Let $(X_j,\omega_j)$, $j=1,2$, have dimensions $n_j\geq 2$ and put $\omega=\omega_1+\omega_2$ and $n=n_1+n_2$. If both $\omega_j$ are balanced hyperbolic, so is $\omega$. If both are Gauduchon hyperbolic, so is $\omega$.\end{The}

\noindent {\it Proof.} In the compact balanced case this is [MP22, Proposition 2.17]; we give the argument in the complete setting and for the Gauduchon case.

\vspace{1ex}

$\bullet$\, {\it Geometry.} By Lemma \ref{Lem:product-geometry}, the product has full bounded geometry, and its universal cover is $\widetilde{X}_1\times\widetilde{X}_2$ with the product metric.

\vspace{1ex}

$\bullet$\, {\it Expansion.} Since $\omega_1,\omega_2$ have even degree, they commute, and $\omega_j^k=0$ for $k>n_j$. In the binomial expansion of $\omega^{n-1}$, only the terms with $k\leq n_1$ and $n-1-k\leq n_2$ survive, that is $k\in\{n_1-1,n_1\}$: $$\omega^{n-1}=a\,\omega_1^{n_1-1}\wedge\omega_2^{n_2}+b\,\omega_1^{n_1}\wedge\omega_2^{n_2-1},\hspace{3ex}a={n-1\choose n_1-1},\ b={n-1\choose n_1}.$$

\vspace{1ex}

$\bullet$\, {\it Both terms are closed (balanced case).} We have $d\omega_j^{n_j-1}=0$ by assumption, and $\omega_j^{n_j}$ is a top-degree form on $X_j$, hence closed. By the Leibniz rule, $d(\omega_1^{n_1-1}\wedge\omega_2^{n_2})=d\omega_1^{n_1-1}\wedge\omega_2^{n_2}=0$, and similarly for the other term. So $\omega$ is balanced.

\vspace{1ex}

$\bullet$\, {\it Bounded primitives.} Choose bounded forms $\beta_j$ with $\widetilde\omega_j^{\,n_j-1}=d\beta_j$. Since $d\widetilde\omega_j^{\,n_j}=0$, the Leibniz rule gives $$d(\beta_1\wedge\widetilde\omega_2^{\,n_2})=\widetilde\omega_1^{\,n_1-1}\wedge\widetilde\omega_2^{\,n_2},\hspace{3ex}d(\widetilde\omega_1^{\,n_1}\wedge\beta_2)=\widetilde\omega_1^{\,n_1}\wedge\widetilde\omega_2^{\,n_2-1},$$ and both primitives are bounded because $|\beta_1\wedge\widetilde\omega_2^{n_2}|=|\beta_1|\,|\widetilde\omega_2^{n_2}|$ and similarly for the second (Lemma \ref{Lem:elementary}(c)). So $\widetilde\omega^{\,n-1}$ is $\widetilde{d}(\mbox{bounded})$.

\vspace{1ex}

$\bullet$\, {\it The Gauduchon case.} Now $\partial\bar\partial\omega_j^{n_j-1}=0$ and $\widetilde\omega_j^{\,n_j-1}=\partial\alpha_j+\bar\partial\beta_j$ with $\alpha_j,\beta_j$ bounded. Because $\omega_j^{n_j}$ has top degree, $\partial\omega_j^{n_j}=\bar\partial\omega_j^{n_j}=0$. Hence $\partial\bar\partial(\omega_1^{n_1-1}\wedge\omega_2^{n_2})=\partial\bar\partial\omega_1^{n_1-1}\wedge\omega_2^{n_2}=0$ and likewise for the other term, so $\partial\bar\partial\omega^{n-1}=0$. For the potentials, $$\widetilde\omega_1^{\,n_1-1}\wedge\widetilde\omega_2^{\,n_2}=\partial(\alpha_1\wedge\widetilde\omega_2^{\,n_2})+\bar\partial(\beta_1\wedge\widetilde\omega_2^{\,n_2}),$$ because $\partial\widetilde\omega_2^{\,n_2}=\bar\partial\widetilde\omega_2^{\,n_2}=0$; and, writing $\widetilde\omega_2^{\,n_2-1}=\partial\alpha_2+\bar\partial\beta_2$, $$\widetilde\omega_1^{\,n_1}\wedge\widetilde\omega_2^{\,n_2-1}=\partial(\widetilde\omega_1^{\,n_1}\wedge\alpha_2)+\bar\partial(\widetilde\omega_1^{\,n_1}\wedge\beta_2),$$ because $\widetilde\omega_1^{\,n_1}$ has even degree and is $\partial$- and $\bar\partial$-closed. All the potentials are bounded as before. \hfill $\qed$

\vspace{2ex}

Before the example we need two facts about the factors.

\begin{Lem}\label{Lem:M_05} Let $\mathcal{M}_{0,5}=\{(x,y)\in\C^2\,;\,x,y\notin\{0,1\},\ x\neq y\}$, the moduli space of five labelled points $0,1,\infty,x,y$ of $\mathbb{P}^1$. Then $\mathcal{M}_{0,5}$ is a smooth complex surface with $\chi(\mathcal{M}_{0,5})=2$, and it carries a K\"ahler metric $\omega_1$ for which $(\mathcal{M}_{0,5},\omega_1)$ is K\"ahler hyperbolic.\end{Lem}

\noindent {\it Proof.} $\bullet$\, {\it Smoothness and Euler characteristic.} $\mathcal{M}_{0,5}$ is the complement of the five lines $x=0$, $x=1$, $y=0$, $y=1$, $x=y$ in $\C^2$, hence a smooth complex surface. Projecting to $x\in\C\setminus\{0,1\}$ (Euler characteristic $-1$), the fibre is $\C\setminus\{0,1,x\}$ (Euler characteristic $-2$), so $\chi(\mathcal{M}_{0,5})=(-1)(-2)=2$.

\vspace{1ex}

$\bullet$\, {\it Uniformisation.} The Teichm\"uller space $\mathcal{T}_{0,5}$ is a simply connected complex surface, and $\mathcal{M}_{0,5}\cong\mathcal{T}_{0,5}/\mbox{PMod}$, where $\mbox{PMod}$ is the pure mapping class group of the five-punctured sphere acting properly discontinuously.

\vspace{1ex}

$\bullet$\, {\it The action is free.} Let a pure mapping class fix the point of $\mathcal{T}_{0,5}$ represented by a marked surface $Y$. It is then induced by a conformal automorphism of $Y$ fixing each of the five punctures. Such an automorphism is a M\"obius transformation of $\mathbb{P}^1$ with five fixed points, hence the identity. So the stabilisers in $\mbox{PMod}$ are trivial, and $\mathcal{T}_{0,5}\to\mathcal{M}_{0,5}$ is the universal covering.

\vspace{1ex}

$\bullet$\, {\it The metric descends.} McMullen's form $\omega_{1/\ell}=\omega_{\rm WP}-i\delta\sum_{\ell_\gamma<\varepsilon}\partial\bar\partial\,\mbox{Log}(\varepsilon/\ell_\gamma)$ is defined from the Weil--Petersson form and the lengths of all closed geodesics; both are invariant under every homeomorphism of the surface, hence under the full mapping class group and in particular under $\mbox{PMod}$. So it defines a K\"ahler form $\omega_1$ on $\mathcal{M}_{0,5}$.

\vspace{1ex}

$\bullet$\, {\it K\"ahler hyperbolicity.} By [McM00], $(\mathcal{T}_{0,5},\widetilde\omega_{1/\ell})$ is complete with bounded curvature, its injectivity radius is bounded below, and $\widetilde\omega_{1/\ell}$ is $d(\mbox{bounded})$; these are (G1), (G3), (G4) and the boundedness of the primitive for $\mathcal{M}_{0,5}$, whose universal cover is $\mathcal{T}_{0,5}$. Finally, $\mathcal{M}_{0,5}$ is a finite cover (of degree $120$) of the orbifold $\mathcal{M}_{0,5}/S_5$, which has finite volume, so (G2) holds. \hfill $\qed$

\begin{Lem}\label{Lem:X_k} For $k\geq 2$, $X_k=\#_k(S^3\times S^3)$ is simply connected with $b_0=b_6=1$, $b_3=2k$ and $b_i=0$ otherwise, so $\chi(X_k)=2-2k$. If $\omega_X$ is a balanced metric on $X_k$ for some complex structure, then $\omega_X^2=d\beta_X$ for a $C^\infty$ $3$-form $\beta_X$, and $X_k$ is degenerate balanced (bounded), hence balanced hyperbolic.\end{Lem}

\noindent {\it Proof.} A connected sum of simply connected closed manifolds is simply connected, and $H_i(M\#N)=H_i(M)\oplus H_i(N)$ for $0<i<6$, with $b_3(S^3\times S^3)=2$ and all other intermediate Betti numbers equal to zero. As $\omega_X$ is balanced, $\omega_X^2$ is a closed $4$-form, and $b_4=0$ makes it exact. On a compact manifold $\beta_X$ is bounded, so Definition \ref{Def:hyperbolic-notions} applies and Proposition \ref{Prop:implications} gives balanced hyperbolicity (compare [MP22, Lemma 2.13]). \hfill $\qed$

\begin{Ex}\label{Ex:product} Let $X_2=\#_2(S^3\times S^3)$ with a balanced metric $\omega_X$ [FLY12], [MP22, Theorem 1.2], and let $X=\mathcal{M}_{0,5}\times X_2$ with $\omega=\omega_1+\omega_X$, so that $n=2+3=5$. Write $\omega_1=d\theta_1$ on $\mathcal{T}_{0,5}$ with $\theta_1$ bounded, and $\omega_X^2=d\beta_X$ on $X_2$ (Lemma \ref{Lem:X_k}). We check the following claims.

\vspace{1ex}

(1)\, {\it $X$ is balanced hyperbolic.} By Lemmas \ref{Lem:M_05} and \ref{Lem:X_k} and Proposition \ref{Prop:implications}, both factors are balanced hyperbolic, so Theorem \ref{The:product-theorem} applies. This is the non-compact analogue of [MP22, Corollary 2.18].

\vspace{1ex}

(2)\, {\it An explicit primitive.} Here $\omega^4=4\,\omega_1\wedge\omega_X^3+6\,\omega_1^2\wedge\omega_X^2$, and $$\widetilde\omega^{\,4}=d\big(4\,\theta_1\wedge\omega_X^3+6\,\omega_1^2\wedge\beta_X\big).$$ Indeed, $d(\theta_1\wedge\omega_X^3)=\omega_1\wedge\omega_X^3$ since $d\theta_1=\omega_1$ and $d\omega_X^3=0$, and $d(\omega_1^2\wedge\beta_X)=\omega_1^2\wedge\omega_X^2$ since $d\omega_1^2=0$ and $d\beta_X=\omega_X^2$. The primitive is bounded as a sum of wedge products of bounded forms.

\vspace{1ex}

(3)\, {\it $X$ carries no K\"ahler metric.} Let $\omega_K$ be a K\"ahler form on $X$ and $i:X_2\to X$, $i(x)=(p,x)$. As $i$ is a holomorphic embedding, $i^\star\omega_K$ is a K\"ahler form on $X_2$. Since $b_2(X_2)=0$, $i^\star\omega_K=d\gamma$ is exact, so $\int_{X_2}(i^\star\omega_K)^3=\int_{X_2}d\big(\gamma\wedge(i^\star\omega_K)^2\big)=0$. But $(i^\star\omega_K)^3$ is a positive multiple of the volume form, a contradiction.

\vspace{1ex}

(4)\, {\it Spectral gap.} We have $\omega^5=10\,\omega_1^2\wedge\omega_X^3$ and $\omega_1^2=d(\theta_1\wedge\omega_1)$, so $$\widetilde\omega^{\,5}=d\big(10\,\theta_1\wedge\omega_1\wedge\omega_X^3\big)$$ has a bounded primitive. Since $\widetilde\omega^5=5!\,dV$, Proposition \ref{Prop:volume-form} gives $\lambda_0(\widetilde{X})>0$. Theorem \ref{The:Cheeger} does not apply here, since $\omega_X$ is not $(\partial+\bar\partial)$-exact on the simply connected compact factor.

\vspace{1ex}

(5)\, {\it Euler characteristic.} $\chi(X)=\chi(\mathcal{M}_{0,5})\,\chi(X_2)=2\cdot(-2)=-4$, which has the sign $(-1)^5$.\end{Ex}

\begin{Prop}\label{Prop:L2-cohomology} Let $X=\mathcal{M}_{0,5}\times X_2$ be as in Example \ref{Ex:product}, with its product metric. The space of $L^2$ harmonic forms on $\widetilde{X}=\mathcal{T}_{0,5}\times X_2$, that is the reduced $L^2$ cohomology, is non-zero in degrees $2$, $5$ and $8$ and zero in all other degrees. In particular, it is not concentrated in the middle degree $n=5$, so Gromov's vanishing theorem does not extend to the class of balanced hyperbolic manifolds.\end{Prop}

\noindent {\it Proof.} $\bullet$\, {\it Product structure.} Since $X_2$ is simply connected, $\widetilde{X}=\mathcal{T}_{0,5}\times X_2$ with the product metric, a complete manifold. On such a manifold the kernel of the Laplacian on forms is the space of $L^2$ harmonic forms, which represents the reduced $L^2$ cohomology.

\vspace{1ex}

$\bullet$\, {\it Separation of variables.} On forms of type $\Lambda^pT^\star\otimes\Lambda^qX_2^\star$ the Laplacian is $A\otimes 1+1\otimes B$ with $A=\Delta_{\mathcal{T}}$ and $B=\Delta_{X_2}$, non-negative self-adjoint operators. If $u$ lies in its kernel, then $0=\langle(A\otimes 1+1\otimes B)u,u\rangle=\|(A^{1/2}\otimes 1)u\|^2+\|(1\otimes B^{1/2})u\|^2$, so $u$ lies in $\ker A\otimes L^2(X_2)$ and in $L^2(\mathcal{T})\otimes\ker B$, that is in $\ker A\otimes\ker B$. The reverse inclusion is clear.

\vspace{1ex}

$\bullet$\, {\it The factors.} By Gromov and McMullen [Gro91], [McM00], $\ker\Delta_{\mathcal{T}}$ vanishes on $p$-forms for $p\neq 2$. For $p=2$ it is non-zero: by the $L^2$ index theorem [CG85], its von Neumann dimension equals $\chi(\mathcal{M}_{0,5})=2>0$. On the compact manifold $X_2$, $\ker\Delta_{X_2}$ on $q$-forms is the space of harmonic forms, of dimension $b_q(X_2)$, which is non-zero exactly for $q=0,3,6$.

\vspace{1ex}

$\bullet$\, {\it Conclusion.} By the two previous steps, the $L^2$ harmonic forms of degree $k=p+q$ are non-zero exactly when $p=2$ and $q\in\{0,3,6\}$, that is $k\in\{2,5,8\}$. \hfill $\qed$

\vspace{2ex}

The von Neumann dimensions in these degrees are $2\cdot 1$, $2\cdot 4$ and $2\cdot 1$, with alternating sum $2-8+2=-4=\chi(X)$, as the $L^2$ index theorem requires.

\section{Entire curves and divisorial hyperbolicity without compactness}\label{section:Brody-div}

We start with two elementary lemmas used in the proof of Theorem \ref{The:Brody-div}.

\begin{Lem}\label{Lem:ODE} (a)\, There is no positive $C^1$ function $A$ on an interval $(r_0,\infty)$ such that $A^2\leq K\,r\,A'$ for some constant $K>0$.

\vspace{1ex}

(b)\, If $G>0$ is $C^1$ on $(r_0,\infty)$ and $G'\geq c_0G$ for some $c_0>0$, then $G(b)\geq G(a)\,e^{c_0(b-a)}$ for all $b>a>r_0$.\end{Lem}

\noindent {\it Proof.} (a) The inequality says $(-1/A)'\geq 1/(Kr)$. Integrating from $r_1$ to $r$ gives $1/A(r_1)>1/A(r_1)-1/A(r)\geq K^{-1}\log(r/r_1)$, which is false for $r$ large. (b) $(e^{-c_0r}G)'\geq 0$. \hfill $\qed$

\begin{Lem}\label{Lem:types} Let $f:\C^{n-1}\to N$ be a holomorphic map and $\alpha$ a $(2n-3)$-form on $N$. Then $f^\star\alpha$ has only components of bidegree $(n-2,n-1)$ and $(n-1,n-2)$. Moreover, $\partial f^\star\alpha=\partial f^\star\alpha^{n-2,n-1}=d f^\star\alpha^{n-2,n-1}$ and $\bar\partial f^\star\alpha=d f^\star\alpha^{n-1,n-2}$.\end{Lem}

\noindent {\it Proof.} On $\C^{n-1}$ a bidegree $(a,b)$ requires $a,b\leq n-1$, and $a+b=2n-3$ leaves only the two pairs stated. Then $\partial f^\star\alpha^{n-1,n-2}$ has bidegree $(n,n-2)$ and vanishes, and $\bar\partial f^\star\alpha^{n-2,n-1}$ has bidegree $(n-2,n)$ and vanishes, which gives the identities. Pull-back by a holomorphic map commutes with $\partial$ and $\bar\partial$. \hfill $\qed$

\vspace{2ex}

\noindent {\it Proof of Theorem \ref{The:Brody-div}(i).} Let $f:\C\to X$ be holomorphic and let $\tilde f:\C\to\widetilde{X}$ be a lift. For a holomorphic map, $\tilde f_\star\partial_y=J\tilde f_\star\partial_x$, and since $\omega(Y,JY)=g(Y,Y)$ we get $$\tilde f^\star\widetilde\omega=|\tilde f'|^2\,dx\wedge dy,\hspace{3ex}|\tilde f'|:=|\tilde f_\star\partial_x|_{\widetilde g}.$$ If $\widetilde\omega=d\theta$ with $\theta$ bounded, put $\eta=\tilde f^\star\theta$. If $\widetilde\omega=\partial\alpha+\bar\partial\beta$ with $\alpha,\beta$ bounded, we argue as in Lemma \ref{Lem:types}. For a $1$-form $\xi$ on $\C$, the form $\partial\xi^{1,0}$ has type $(2,0)$ and $\bar\partial\xi^{0,1}$ has type $(0,2)$, so both vanish, and therefore $\partial\xi=d\xi^{0,1}$ and $\bar\partial\xi=d\xi^{1,0}$. As pull-back by a holomorphic map commutes with $\partial$ and $\bar\partial$, we get $\tilde f^\star\widetilde\omega=d\eta$ with $\eta=\tilde f^\star\alpha^{0,1}+\tilde f^\star\beta^{1,0}$. In both cases, $|\eta(v)|\leq B\,|\tilde f'|$ for every unit vector $v$, where $B$ bounds the sup-norms of $\theta$, or of $\alpha$ and $\beta$.

Suppose $f$ is non-constant. Then $A(r)=\int_{D_r}|\tilde f'|^2dx\,dy$ is positive for $r>r_0$, it is $C^1$, and $A'(r)=\int_0^{2\pi}|\tilde f'(re^{it})|^2\,r\,dt$. By Stokes' theorem and the Cauchy--Schwarz inequality, $$A(r)=\int_{\partial D_r}\eta\leq B\int_0^{2\pi}|\tilde f'(re^{it})|\,r\,dt\leq B\sqrt{2\pi r\,A'(r)}.$$ Thus $A^2\leq 2\pi B^2\,r\,A'$, which contradicts Lemma \ref{Lem:ODE}(a). \hfill $\qed$

\vspace{2ex}

For (ii) we recall the notion of subexponential growth of [MP22, Definition 2.3], which is stated there for compact manifolds but makes sense verbatim for any Hermitian manifold $(X,\omega)$. Let $f:\C^{n-1}\to X$ be a holomorphic map that is non-degenerate at some point. Put $\omega_{n-1}=\omega^{n-1}/(n-1)!$, let $B_t$ and $S_t$ be the ball and the sphere of radius $t$ in $\C^{n-1}$, $\tau=|z|^2$, $$\mbox{Vol}_{\omega,f}(B_t)=\int_{B_t}f^\star\omega_{n-1},\hspace{3ex}F(b)=\int_0^b\mbox{Vol}_{\omega,f}(B_t)\,dt,$$ and let $d\sigma_{\omega,f,t}$ be the area measure induced by $f^\star\omega$ on $S_t$ (off the proper analytic set where $f$ is degenerate).

\begin{Def}\label{Def:subexp-growth} The map $f$ has {\bf subexponential growth} with respect to $\omega$ if:

\vspace{1ex}

(i)\, there exist $C_1,r_0>0$ such that $\displaystyle\int_{S_t}\|d\tau\|_{f^\star\omega}\,d\sigma_{\omega,f,t}\leq C_1\,t\,\mbox{Vol}_{\omega,f}(B_t)$ for all $t>r_0$;

\vspace{1ex}

(ii)\, for every $C>0$, $\displaystyle\limsup_{b\to+\infty}\Big(\frac{b}{C}-\log F(b)\Big)=+\infty$.

\vspace{1ex}

\noindent The manifold $(X,\omega)$ is said to be {\bf divisorially hyperbolic with respect to $\omega$} if there is no holomorphic map $f:\C^{n-1}\to X$, non-degenerate at some point, of subexponential growth with respect to $\omega$.\end{Def}

If $\omega_1,\omega_2$ are uniformly equivalent, the two notions of subexponential growth coincide, since by Lemma \ref{Lem:elementary}(b) every quantity involved changes by a bounded factor. For compact $X$ this is [MP22, Lemma 2.6].

\vspace{2ex}

\noindent {\it Proof of Theorem \ref{The:Brody-div}(ii).} The argument follows [MP22, Theorem 2.8], where compactness is used only to make the notion independent of the metric. Let $\tilde f:\C^{n-1}\to\widetilde{X}$ be a lift of $f$.

\vspace{1ex}

\noindent {\bf Claim.} {\it There is a bounded $(2n-3)$-form $\widehat\Gamma$ on $\widetilde{X}$ such that $\tilde f^\star\widetilde\omega^{\,n-1}=d(\tilde f^\star\widehat\Gamma)$.}

\vspace{1ex}

\noindent {\it Proof of the Claim.} In the $\widetilde{d}(\mbox{bounded})$ case, take $\widehat\Gamma=\Gamma$, where $\widetilde\omega^{n-1}=d\Gamma$. If $\widetilde\omega^{n-1}=\partial\alpha+\bar\partial\beta$, Lemma \ref{Lem:types} gives $\tilde f^\star\widetilde\omega^{\,n-1}=d\tilde f^\star(\alpha^{n-2,n-1}+\beta^{n-1,n-2})$, and we take $\widehat\Gamma=\alpha^{n-2,n-1}+\beta^{n-1,n-2}$, which is bounded by Lemma \ref{Lem:elementary}(a). \hfill $\qed$

\vspace{1ex}

\noindent $\bullet$\, {\it Non-degenerate maps from $\mathbb{P}^{n-1}$.} Since $\mathbb{P}^{n-1}$ is simply connected, such a map lifts to $\widetilde{X}$ and the Claim holds there, so $\int_{\mathbb{P}^{n-1}}f^\star\omega^{n-1}=\int d(\cdots)=0$. This contradicts the positivity of $f^\star\omega^{n-1}$ where $f$ is non-degenerate.

\vspace{1ex}

\noindent $\bullet$\, {\it Growth.} Suppose $f:\C^{n-1}\to X$ is non-degenerate at some point and satisfies (i) and (ii). Take $r_0$ larger than the norm of that point, so that $V(t):=\mbox{Vol}_{\omega,f}(B_t)>0$ for $t>r_0$, and write $A(t)=\mbox{area}_{f^\star\omega}(S_t)$. By construction $|v_j|_{f^\star\omega}=|\tilde f_\star v_j|_{\widetilde\omega}$, so for tangent vectors $v_1,\dots,v_{2n-3}$ to $S_t$ we have $|\tilde f^\star\widehat\Gamma(v_1,\dots,v_{2n-3})|\leq\|\widehat\Gamma\|_\infty\prod|v_j|_{f^\star\omega}$. Hence, by Stokes' theorem (the set where $f$ is degenerate meets $S_t$ in a set of measure zero), $$V(t)=\frac{1}{(n-1)!}\int_{S_t}\tilde f^\star\widehat\Gamma\leq C\,A(t).\hspace{3ex}(1)$$ Next, we use the coarea formula for $\tau$ on the set where $g=f^\star\omega$ is Riemannian; there $f^\star\omega_{n-1}$ is the Riemannian volume form of $g$, and $d\tau=2t\,dt$, so $$V(r)=\int_0^r\Big(\int_{S_t}\frac{d\sigma_{\omega,f,t}}{\|d\tau\|_g}\Big)\,2t\,dt.$$ The Cauchy--Schwarz inequality gives $A(t)^2\leq\int_{S_t}\|d\tau\|_g^{-1}d\sigma\cdot\int_{S_t}\|d\tau\|_g\,d\sigma$. Combining with (1) and (i), for $t>r_0$, $$\int_{S_t}\frac{d\sigma_{\omega,f,t}}{\|d\tau\|_g}\geq\frac{A(t)^2}{\int_{S_t}\|d\tau\|_g\,d\sigma_{\omega,f,t}}\geq\frac{V(t)^2/C^2}{C_1\,t\,V(t)}=\frac{V(t)}{C^2C_1\,t}.$$ Therefore, with $c_0=2/(C^2C_1)$ and $G(r)=\int_{r_0}^rV(t)\,dt>0$, we get $V(r)\geq c_0G(r)$, that is $G'\geq c_0G$. By Lemma \ref{Lem:ODE}(b), $\log G(b)\geq\log G(a)+c_0(b-a)$ for $b>a>r_0$. Since $F\geq G$, the quantity $c_0b-\log F(b)$ is bounded above as $b\to\infty$, which contradicts (ii) with $C=1/c_0$. \hfill $\qed$

\vspace{2ex}

\begin{Prop}\label{Prop:Brody-div} If $X$ is Brody hyperbolic, no holomorphic map $\C^{n-1}\to X$ is non-degenerate at any point. In particular, $(X,\omega)$ is divisorially hyperbolic with respect to every Hermitian metric $\omega$.\end{Prop}

\noindent {\it Proof.} Let $f:\C^{n-1}\to X$ be non-degenerate at $x$. Since $n\geq 2$, $d_xf$ is injective on a complex line $L$ through $x$, and $f|_L:L\to X$ is a non-constant entire curve. \hfill $\qed$

\vspace{2ex}

For non-compact $X$, Brody hyperbolicity does not imply Kobayashi hyperbolicity, and Brody's theorem [Bro78] identifying the two needs compactness; so part (i) of Theorem \ref{The:Brody-div} gives Brody hyperbolicity only.

\begin{Question}\label{Question:Kobayashi} Is a complete K\"ahler hyperbolic manifold with full bounded geometry Kobayashi hyperbolic? The proof of part (i) bounds the area of $f(D_{1/2})$ for every holomorphic disc $f$ by a universal constant. A monotonicity formula and Cauchy estimates in uniform holomorphic charts would then bound the Kobayashi metric from below, but (G3) and (G4) do not obviously provide such charts.\end{Question}

\section{A non-compact surface}\label{section:surface}

For compact surfaces, divisorial and Brody hyperbolicity coincide by [Mar26, Proposition 3.1]: starting from a Brody curve $g$ with $|g'|\leq 1$, one composes with $z\mapsto z^k$ or with $\exp$ to obtain a curve of subexponential growth. Compactness enters through Brody's reparametrisation lemma, and the equivalence fails without it. We need two elementary lemmas: a computation of condition (i) in dimension two, and a quantitative form of Bloch's theorem.

\begin{Lem}\label{Lem:dim-two} Let $n=2$ and let $f:\C\to X$ be a non-constant holomorphic map, with $A_f(t)=\int_{D_t}|f'|^2\,dx\,dy$ and $|f'|=|f_\star\partial_x|_\omega$. Then, for every $t>0$, $$\int_{S_t}\|d\tau\|_{f^\star\omega}\,d\sigma_{\omega,f,t}=4\pi t^2.$$ Consequently, condition (i) of Definition \ref{Def:subexp-growth} holds if and only if $A_f(t)\geq c\,t$ for some $c>0$ and all large $t$.\end{Lem}

\noindent {\it Proof.} As in the proof of Theorem \ref{The:Brody-div}(i), $f^\star\omega=\lambda^2|dz|^2$ with $\lambda=|f'|$. Off the discrete set where $\lambda=0$, the metric $g=f^\star\omega$ is conformal to the Euclidean metric, so $\|d\tau\|_g=\|d\tau\|_{\rm eucl}/\lambda=2t/\lambda$, and the induced length element on $S_t$ is $d\sigma_{\omega,f,t}=\lambda\,ds$ with $ds$ the Euclidean arc length. The product is $2t\,ds$ and its integral over $S_t$ is $4\pi t^2$. Since $n=2$ we have $\omega_{n-1}=\omega$, so $\mbox{Vol}_{\omega,f}(B_t)=A_f(t)$, and (i) reads $4\pi t^2\leq C_1\,t\,A_f(t)$, that is $A_f(t)\geq(4\pi/C_1)\,t$. \hfill $\qed$

\begin{Lem}\label{Lem:Bloch} Let $h:\C\to\C$ be a non-constant entire function and $\rho(w)=e^{|w|^2}$. There exist constants $c,c'>0$ and $t_1$ such that, for $t\geq t_1$, $$\int_{D_t}\rho(h(z))^2\,|h'(z)|^2\,dx\,dy\geq c\,e^{c't^2}.$$\end{Lem}

\noindent {\it Proof.} Choose $z_0$ with $s=|h'(z_0)|>0$. Bloch's theorem asserts that there is an absolute constant $b>0$ such that, if $\phi$ is holomorphic on the unit disc with $|\phi'(0)|=1$, then $\phi(\mathbb{D})$ contains a disc of radius $b$. Applying it to $\phi(\zeta)=h(z_0+r\zeta)/(rs)$ shows that $h(D_r(z_0))$ contains a Euclidean disc $\Delta$ of radius $b\,s\,r$. For $t>2|z_0|$ we have $D_{t/2}(z_0)\subset D_t$, so $h(D_t)\supset\Delta$ with $R:=bst/2$ the radius of $\Delta$. The integral in the lemma is the $\rho^2$-area of $h(D_t)$ counted with multiplicity, hence at least the $\rho^2$-area of $\Delta$. A disc of radius $R$ meets the disc $\{|w|<R/2\}$ in an area at most $\pi R^2/4$, so at least $3/4$ of the area of $\Delta$ lies in $\{|w|\geq R/2\}$, where $\rho^2\geq e^{R^2/2}$. Thus the integral is at least $\frac{3\pi}{4}R^2e^{R^2/2}\geq c\,e^{c't^2}$ with $c'=b^2s^2/8$. \hfill $\qed$

\vspace{2ex}

\noindent {\it Proof of Theorem \ref{The:surface}.} Let $X=\C\times\mathbb{D}$ with $\omega=\rho(z)^2|dz|^2+\omega_{\mathbb{D}}$, $\rho(z)=e^{|z|^2}$.

\vspace{1ex}

$\bullet$\, {\it Geometry of $(X,\omega)$.} It is a product of two K\"ahler metrics, hence K\"ahler. The metric $\rho^2|dz|^2$ on $\C$ is complete since $\rho\geq 1$, and its curvature is $-\Delta\log\rho/\rho^2=-4e^{-2|z|^2}\in[-4,0)$ because $\Delta|z|^2=4$. The Poincar\'e metric is complete of curvature $-1$. By the curvature computation in the proof of Lemma \ref{Lem:product-geometry}, $\omega$ is complete with sectional curvature in $[-4,0]$; it is simply connected, so by the Cartan--Hadamard theorem $X$ is a Hadamard manifold with infinite injectivity radius. Its volume is infinite.

\vspace{1ex}

$\bullet$\, {\it $X$ is not Brody hyperbolic}, since $z\mapsto(z,0)$ is an entire curve.

\vspace{1ex}

$\bullet$\, {\it $X$ is divisorially hyperbolic with respect to $\omega$.} Let $f=(h,c):\C\to X$ be non-constant. A holomorphic map $\C\to\mathbb{D}$ is constant by Liouville's theorem, so $c$ is constant and $h$ is a non-constant entire function. By Lemma \ref{Lem:Bloch}, $A_f(t)=\int_{D_t}\rho(h)^2|h'|^2\geq c\,e^{c't^2}$ for large $t$. As $A_f$ is non-decreasing, $F(b)\geq A_f(b-1)$, hence $\log F(b)\geq c'(b-1)^2-O(1)$ and $b/C-\log F(b)\to-\infty$ for every $C>0$. So condition (ii) fails for every non-constant $f$.

\vspace{1ex}

$\bullet$\, {\it The metric $\omega'=|dz|^2+\omega_{\mathbb{D}}$.} For $f(z)=(z,0)$ we have $A_f(t)=\pi t^2$. By Lemma \ref{Lem:dim-two}, condition (i) amounts to $A_f(t)\geq c\,t$ for large $t$, which holds, and $F(b)=\pi b^3/3$ gives (ii). Hence $(X,\omega')$ is not divisorially hyperbolic. This is consistent with Definition \ref{Def:subexp-growth}, since $\omega$ and $\omega'$ are not uniformly equivalent. \hfill $\qed$

\begin{Question}\label{Question:surfaces} Are Brody and divisorial hyperbolicity equivalent for complex surfaces with full bounded geometry? The example above has infinite volume, and the compact proof uses Brody's lemma, which needs compactness.\end{Question}

\section{Functorial pseudo-distances}\label{section:pseudo-distances}

The Kobayashi pseudo-distance $d_K$ [Kob70] is functorial: holomorphic maps do not increase it, and $X$ is Kobayashi hyperbolic exactly when $d_K$ is a distance. This is what makes Kobayashi hyperbolicity easy to handle; it is inherited by complex submanifolds, it gives Brody's theorem on compact manifolds, and it underlies the Picard-type theorems. Divisorial hyperbolicity has no such formulation. It is defined by a growth condition on maps $\C^{n-1}\to X$, relative to a metric, and it is natural to ask for an analogue of $d_K$, built from such maps, which would be non-degenerate exactly on divisorially hyperbolic manifolds. Theorem \ref{The:pseudo-distance} says that no functorial continuous construction can do this. The reason is simple: a functorial pseudo-distance vanishes along every entire curve, so it can only detect Brody-type properties, whereas divisorial hyperbolicity allows entire curves. This also explains why the growth conditions of Definition \ref{Def:subexp-growth}, and the dependence on the metric, cannot be avoided.

\vspace{2ex}

\noindent {\it Proof of Theorem \ref{The:pseudo-distance}.} $\bullet$\, {\it $d_\C\equiv 0$.} Let $x\neq y$ and $x'\neq y'$ in $\C$. The affine map $\varphi(z)=x'+(y'-x')(z-x)/(y-x)$ is holomorphic with $\varphi(x)=x'$ and $\varphi(y)=y'$, so $d_\C(x',y')\leq d_\C(x,y)$. Exchanging the two pairs gives equality, so $d_\C$ takes a constant value $c$ on distinct pairs. By continuity at the diagonal, where $d_\C$ vanishes, $c=0$.

\vspace{1ex}

$\bullet$\, {\it Consequences.} If $f:\C\to X$ is holomorphic, $d_X(f(a),f(b))\leq d_\C(a,b)=0$. If $f$ is non-constant, it takes two distinct values at distance $0$, so $d_X$ is not a distance. Hence a distance $d_X$ forces $X$ to be Brody hyperbolic.

\vspace{1ex}

$\bullet$\, {\it A compact divisorially hyperbolic manifold that is not Brody hyperbolic.} Let $X_1$ be a compact Brody hyperbolic manifold of dimension $n_1\geq 2$ (for instance a smooth ample surface in a simple abelian threefold, by Example \ref{Ex:not-SKT-hyperbolic}(a)) and $X=X_1\times\mathbb{P}^m$ with $m\geq 1$, so that $n=n_1+m$. Let $f=(f_1,f_2):\C^{n-1}\to X$ be holomorphic. The restriction of $f_1$ to any complex line is an entire curve in $X_1$, hence constant; since any two points of $\C^{n-1}$ lie on a line, $f_1$ is constant. Thus $\mbox{rank}\,df\leq\mbox{rank}\,df_2\leq m=n-n_1\leq n-2<n-1$, and no holomorphic map $\C^{n-1}\to X$ is non-degenerate at any point. So $X$ is divisorially hyperbolic (compare [Mar26, Lemma 6.7]). But $\mathbb{P}^m\subset X$ gives non-constant entire curves, so $d_X$ is not a distance. Another example is $SL(2,\C)/\Gamma$ for a cocompact lattice $\Gamma$, which is degenerate balanced [Yac98], hence divisorially hyperbolic by [MP22, Theorem 1.2], and carries entire curves. \hfill $\qed$

\vspace{2ex}

We conclude with two statements that clarify the position of divisorial hyperbolicity among functorial notions. The first says that chains of maps from balls give nothing new.

\begin{Prop}\label{Prop:balls} For $x,y\in X$, let $d^{(n-1)}_X(x,y)$ be the infimum of $\sum_{i=1}^kd_B(a_i,b_i)$ over all chains of holomorphic maps $f_i:B^{n-1}\to X$ and points $a_i,b_i\in B^{n-1}$ with $f_1(a_1)=x$, $f_i(b_i)=f_{i+1}(a_{i+1})$ and $f_k(b_k)=y$, where $d_B$ is the Kobayashi distance of the ball. Then $d^{(n-1)}_X=d^K_X$.\end{Prop}

\noindent {\it Proof.} Each $f_i$ decreases Kobayashi distances, so $d^K_X(f_i(a_i),f_i(b_i))\leq d_B(a_i,b_i)$, and the triangle inequality gives $d^K_X(x,y)\leq d^{(n-1)}_X(x,y)$. Conversely, take a chain of discs $g_i:\mathbb{D}\to X$ with points $a_i,b_i\in\mathbb{D}$ computing $d^K_X(x,y)$ up to $\varepsilon$, and put $f_i=g_i\circ\mbox{pr}_1$ and $\iota(\zeta)=(\zeta,0)$. Since $\mbox{pr}_1\circ\iota=\mbox{id}$, both $\iota$ and $\mbox{pr}_1$ are distance-decreasing and $d_B(\iota a,\iota b)=d_{\mathbb{D}}(a,b)$. The chain $(f_i,\iota a_i,\iota b_i)$ therefore has the same length, so $d^{(n-1)}_X\leq d^K_X$. \hfill $\qed$

\vspace{2ex}

Restricting to non-degenerate maps, or to maps of subexponential growth, changes the infimum, but these classes are not stable under composition, so functoriality is lost. The second statement shows where the functorial condition ``no non-degenerate map $\C^{n-1}\to X$'' sits.

\begin{Prop}\label{Prop:in-between} Call $X$ {\bf $(n-1)$-Brody hyperbolic} if no holomorphic map $\C^{n-1}\to X$ is non-degenerate at some point. Then $$\mbox{Brody hyperbolic}\ \Longrightarrow\ (n-1)\mbox{-Brody hyperbolic}\ \Longrightarrow\ \mbox{divisorially hyperbolic},$$ and neither implication can be reversed: $X_1\times\mathbb{P}^m$ ($n_1\geq 2$, $m\geq 1$) is $(n-1)$-Brody hyperbolic but not Brody hyperbolic, and $SL(2,\C)/\Gamma$ ($\Gamma$ a cocompact lattice) is divisorially hyperbolic but not $(n-1)$-Brody hyperbolic.\end{Prop}

\noindent {\it Proof.} The first implication is Proposition \ref{Prop:Brody-div}; the second holds by definition, since a map of subexponential growth, non-degenerate at some point, is in particular a non-degenerate map. For $X_1\times\mathbb{P}^m$, the rank computation in the proof of Theorem \ref{The:pseudo-distance} shows that no map $\C^{n-1}\to X$ is non-degenerate, while $\mathbb{P}^m\subset X$ gives entire curves. The manifold $SL(2,\C)/\Gamma$ is divisorially hyperbolic by [MP22, Theorem 1.2], but the map $f(z_1,z_2)=\left(\begin{smallmatrix}e^{z_1}&z_2\\0&e^{-z_1}\end{smallmatrix}\right)\Gamma$ from $\C^2$ is non-degenerate at every point, because $\C^2\to SL(2,\C)$ is, and $SL(2,\C)\to SL(2,\C)/\Gamma$ is a local biholomorphism [MP22, \S2]. \hfill $\qed$

\section{Complements}\label{section:complements}

We say that a $(p,\,p)$-form $\Omega$ on $(X,\omega)$ is {\it uniformly positive} if it is bounded and there exists $c>0$ such that $\Omega|_P\geq c\,\omega^p|_P$ for every complex $p$-dimensional subspace $P$ of a tangent space (as top-degree forms on the oriented space $P$).

\begin{Prop}\label{Prop:isoperimetric} Let $\Omega$ be uniformly positive with $\widetilde\Omega=\partial\alpha+\bar\partial\beta$, $\alpha,\beta$ bounded. There exists $C>0$ such that, for every closed complex submanifold $Z'\subset\widetilde{X}$ of dimension $p$ and every relatively compact open set $F\subset Z'$ with smooth boundary, $\mbox{vol}_{2p}(F)\leq C\,\mbox{vol}_{2p-1}(\partial F)$.\end{Prop}

\noindent {\it Proof.} Since $\widetilde\Omega$ has type $(p,\,p)$, it equals the $(p,\,p)$-part of $\partial\alpha+\bar\partial\beta$, namely $\partial\alpha^{p-1,p}+\bar\partial\beta^{p,p-1}$. On $Z'$, which has complex dimension $p$, the form $\bar\partial\alpha^{p-1,p}$ has type $(p-1,p+1)$ and $\partial\beta^{p,p-1}$ has type $(p+1,p-1)$, so both vanish; hence $\partial\alpha^{p-1,p}|_{Z'}=d\alpha^{p-1,p}|_{Z'}$, and likewise for $\beta$. Put $\eta=\alpha^{p-1,p}+\beta^{p,p-1}$, which is bounded by Lemma \ref{Lem:elementary}(a); then $\widetilde\Omega|_{Z'}=d(\eta|_{Z'})$. By uniform positivity and Wirtinger's identity, $\widetilde\Omega|_{Z'}\geq c\,\widetilde\omega^p|_{Z'}=c\,p!\,dV_{Z'}$, and by Stokes' theorem $$c\,p!\,\mbox{vol}(F)\leq\int_F\widetilde\Omega=\int_{\partial F}\eta\leq\|\eta\|_\infty\,\mbox{vol}_{2p-1}(\partial F).\hspace{3ex}\hfill\qed$$

\begin{Cor}\label{Cor:obstruction} Let $X$ be a compact complex manifold and $Z\subset X$ a compact connected complex submanifold of dimension $p$, $1\leq p\leq n-1$, such that the image $\Gamma'$ of $\pi_1(Z)\to\pi_1(X)$ is amenable. Then $X$ is not $p$-SKT hyperbolic. In particular, for $p=n-1$, $X$ is neither Gauduchon hyperbolic nor balanced hyperbolic.\end{Cor}

\noindent {\it Proof.} Suppose $X$ is $p$-SKT hyperbolic.

\vspace{1ex}

$\bullet$\, {\it The covering.} Let $Z'$ be a connected component of the preimage of $Z$ in $\widetilde{X}$. It is a closed complex submanifold of dimension $p$, and $Z'\to Z$ is a normal covering whose deck group is the image $\Gamma'$ of $\pi_1(Z)\to\pi_1(X)$. By Proposition \ref{Prop:isoperimetric}, $\mbox{vol}(F)\leq C\,\mbox{area}(\partial F)$ for every relatively compact open set $F\subset Z'$ with smooth boundary.

\vspace{1ex}

$\bullet$\, {\it $\Gamma'$ finite.} Then $Z'$ is compact without boundary, and $F=Z'$ gives $\mbox{vol}(Z')\leq 0$, which is absurd.

\vspace{1ex}

$\bullet$\, {\it $\Gamma'$ infinite and amenable.} By Brooks' theorem for normal coverings of a compact manifold [Bro81], [Bro85], $\lambda_0(Z')=0$. If the Cheeger constant of $Z'$ were positive, Cheeger's inequality would give $\lambda_0(Z')\geq h^2/4>0$; hence $h(Z')=0$, so there are domains $F$ with $\mbox{area}(\partial F)/\mbox{vol}(F)$ arbitrarily small, contradicting the first step.

\vspace{1ex}

\noindent The assertion for $p=n-1$ follows since Gauduchon and balanced hyperbolicity imply $(n-1)$-SKT hyperbolicity with $\Omega=\omega^{n-1}$ (Proposition \ref{Prop:implications}). \hfill $\qed$

\begin{Prop}\label{Prop:Poincare} Let $\widetilde\omega^{n-1}=d\beta$ with $\beta$ bounded on $\widetilde{X}$, and let $d^c=i(\bar\partial-\partial)$.

\vspace{1ex}

(a)\, There exists $C>0$ such that $\|df\|_{L^2}\leq C\,\|dd^cf\|_{L^2}$ for all real $f\in C^\infty_c(\widetilde{X})$.

\vspace{1ex}

(b)\, Write $d\widetilde\omega\wedge\beta=u\,\widetilde\omega^n$. If $\sup|u|\leq\varepsilon<1$, then $\lambda_0(\widetilde{X})>0$.\end{Prop}

\noindent {\it Proof.} (a) For real $f$ one has $df\wedge d^cf=2i\,\partial f\wedge\bar\partial f$, so $df\wedge d^cf\wedge\widetilde\omega^{n-1}=c_n|df|^2dV$ with $c_n>0$ depending only on $n$. Since $d(df\wedge d^cf)=-df\wedge dd^cf$ and $f$ has compact support, $$0=\int d\big(df\wedge d^cf\wedge\beta\big)=-\int df\wedge dd^cf\wedge\beta+\int df\wedge d^cf\wedge d\beta.$$ The last integral equals $c_n\|df\|_{L^2}^2$, and the first is bounded by $C'\|\beta\|_\infty\|df\|_{L^2}\|dd^cf\|_{L^2}$. Divide by $\|df\|_{L^2}$.

\vspace{1ex}

(b) From $d(\widetilde\omega\wedge\beta)=d\widetilde\omega\wedge\beta+\widetilde\omega^n$ we get $\widetilde\omega^n=d(\widetilde\omega\wedge\beta)-u\,\widetilde\omega^n$. Multiplying by $f^2$ and integrating, and using $\widetilde\omega^n=n!\,dV$, $$\int f^2\,\widetilde\omega^n=-2\int f\,df\wedge\widetilde\omega\wedge\beta-\int f^2u\,\widetilde\omega^n,$$ so $(1-\varepsilon)\,n!\,\|f\|_{L^2}^2\leq 2c_n'\|\beta\|_\infty\|f\|_{L^2}\|df\|_{L^2}$, whence $\|df\|_{L^2}\geq\frac{(1-\varepsilon)\,n!}{2c_n'\|\beta\|_\infty}\|f\|_{L^2}$. \hfill $\qed$

\vspace{6ex}

\noindent {\bf References.} \\

\noindent [Bro78]\, R. Brody --- {\it Compact Manifolds and Hyperbolicity} --- Trans. Amer. Math. Soc. {\bf 235} (1978), 213-219.

\vspace{1ex}

\noindent [Bro81]\, R. Brooks --- {\it The Fundamental Group and the Spectrum of the Laplacian} --- Comment. Math. Helv. {\bf 56} (1981), 581-598.

\vspace{1ex}

\noindent [Bro85]\, R. Brooks --- {\it The Bottom of the Spectrum of a Riemannian Covering} --- J. Reine Angew. Math. {\bf 357} (1985), 101-114.

\vspace{1ex}

\noindent [CG85]\, J. Cheeger, M. Gromov --- {\it On the Characteristic Numbers of Complete Manifolds of Bounded Curvature and Finite Volume} --- in Differential Geometry and Complex Analysis, Springer-Verlag (1985), 115-154.

\vspace{1ex}

\noindent [Che70]\, J. Cheeger --- {\it A Lower Bound for the Smallest Eigenvalue of the Laplacian} --- in Problems in Analysis, Princeton Univ. Press (1970), 195-199.

\vspace{1ex}

\noindent [CY17]\, B.-L. Chen, X. Yang --- {\it Compact K\"ahler Manifolds Homotopic to Negatively Curved Riemannian Manifolds} --- Math. Ann. {\bf 370} (2017), DOI 10.1007/s00208-017-1521-7

\vspace{1ex}

\noindent [FLY12]\, J. Fu, J. Li, S.-T. Yau --- {\it Balanced Metrics on Non-K\"ahler Calabi-Yau Threefolds} --- J. Differential Geom. {\bf 90} (2012), 81-129.

\vspace{1ex}

\noindent [Fri89]\, R. Friedman --- {\it On Threefolds with Trivial Canonical Bundle} --- in Complex Geometry and Lie Theory (Sundance, UT, 1989), Proc. Sympos. Pure Math., vol. 53, Amer. Math. Soc., Providence, RI, 1991, 103-134.

\vspace{1ex}

\noindent [Gre75]\, M. Green --- {\it Some Picard Theorems for Holomorphic Maps to Algebraic Varieties} --- Amer. J. Math. {\bf 97} (1975), 43-75.

\vspace{1ex}

\noindent [Gro91]\, M. Gromov --- {\it K\"ahler Hyperbolicity and $L^2$ Hodge Theory} --- J. Diff. Geom. {\bf 33} (1991), 263-292.

\vspace{1ex}

\noindent [Khe26]\, A. Khelifati --- {\it Deformations of Non-K\"ahler Hyperbolicity Notions and Modifications of Degenerate Balanced Manifolds} --- arXiv e-print 2512.19284.

\vspace{1ex}

\noindent [Kob70]\, S. Kobayashi --- {\it Hyperbolic Manifolds and Holomorphic Mappings} --- Marcel Dekker, New York (1970).

\vspace{1ex}

\noindent [LT93]\, P. Lu, G. Tian --- {\it The Complex Structures on Connected Sums of $S^3\times S^3$} --- in Manifolds and Geometry (Pisa, 1993), Sympos. Math., XXXVI, Cambridge Univ. Press, Cambridge, 1996, 284-293.

\vspace{1ex}

\noindent [Mar24]\, S. Marouani --- {\it SKT Hyperbolic and Gauduchon Hyperbolic Compact Complex Manifolds} --- Bull. Soc. Math. France {\bf 152}, no. 3 (2024), 519-545.

\vspace{1ex}

\noindent [Mar26]\, S. Marouani --- {\it Non-K\"ahler Hyperbolicity: Surfaces and Codimension One} --- arXiv e-print math.CV 2611.07098.

\vspace{1ex}

\noindent [McM00]\, C. T. McMullen --- {\it The Moduli Space of Riemann Surfaces is K\"ahler Hyperbolic} --- Ann. of Math. {\bf 151} (2000), 327-357.

\vspace{1ex}

\noindent [MP22]\, S. Marouani, D. Popovici --- {\it Balanced Hyperbolic and Divisorially Hyperbolic Compact Complex Manifolds} --- Math. Res. Lett. {\bf 30}, no. 6 (2023), 1813-1855.

\vspace{1ex}

\noindent [Och77]\, T. Ochiai --- {\it On Holomorphic Curves in Algebraic Varieties with Ample Irregularity} --- Invent. Math. {\bf 43} (1977), 83-96.

\vspace{1ex}

\noindent [Yac98]\, A. Yachou --- {\it Sur les vari\'et\'es semi-k\"ahl\'eriennes} --- PhD Thesis, University of Lille.

\vspace{4ex}

\noindent Institut Sup\'erieur des \'Etudes Appliqu\'ees en Humanit\'es de Sbe\"\i tla,

\noindent Universit\'e de Kairouan, Tunisie

\noindent Email: marouanisamir20@gmail.com

\end{document}